\documentclass[letterpaper, 10 pt, conference]{ieeeconf}  

\usepackage{url}  
\usepackage{graphicx}  
\usepackage{amsmath,amssymb,amsfonts,mathrsfs,bm} 
\usepackage{caption}
\usepackage{cleveref}
\usepackage{array}
\usepackage{scalerel}
\usepackage{algorithm}
\usepackage{algpseudocode}
\usepackage{booktabs}
\usepackage{epstopdf}
\usepackage{multirow}
\usepackage{rotating}
\usepackage{algpseudocode}
\usepackage{cite}
\usepackage{epsfig}
\usepackage{tabularx}
\usepackage[table]{xcolor}
\usepackage{epstopdf}
\usepackage{accents}
\usepackage{mathtools}

\newcommand{\msh}{\!\!\;}

\ifx\theorem\undefined

\newtheorem{theorem}{Theorem}
\newenvironment{theorem*}{\par\noindent{\bf Theorem\ }}{\hfill\\[2mm]}

\newenvironment{corollary*}{\par\noindent{\bf Corollary\ }}{\hfill\\[2mm]}
\newtheorem{lemma}{Lemma}
\newtheorem{definition}{Definition}

\fi
\newcommand{\ubar}[1]{\underaccent{\bar}{#1}}

\newcommand{\tr}{\mathrm{tr}}

\newcommand{\abf}{\mathbf{a}}
\newcommand{\Abf}{\mathbf{A}}
\newcommand{\Acal}{\mathcal{A}}

\newcommand{\bbf}{\mathbf{b}}
\newcommand{\Bbf}{\mathbf{B}}

\newcommand{\cbf}{\mathbf{c}}
\newcommand{\Ccal}{\mathcal{C}}

\newcommand{\ebf}{\mathbf{e}}

\newcommand{\Fbf}{\mathbf{F}}

\newcommand{\Fcal}{\mathcal{F}}

\newcommand{\Hbf}{\mathbf{H}}

\newcommand{\Ibf}{\mathbf{I}}

\newcommand{\Kbf}{\mathbf{K}}

\newcommand{\Lbf}{\mathbf{L}}

\newcommand{\Lcal}{\mathcal{L}}

\newcommand{\Mcal}{\mathcal{M}}

\newcommand{\Ncal}{\mathcal{N}}

\newcommand{\Rbb}{\mathbb{R}}

\newcommand{\Sbb}{\mathbb{S}}

\newcommand{\ubf}{\mathbf{u}}

\newcommand{\xbf}{\mathbf{x}}
\newcommand{\Xbf}{\mathbf{X}}

\newcommand{\Lambdabf}{\mathbf{\Lambda}}

\newcommand{\sm}[2]{\scaleto{#1\mathstrut}{#2pt}}
\definecolor{pinkF}{rgb}{0.858, 0.188, 0.478}
\newif\ifcomment

\crefrangelabelformat{equation}{(#3#1#4)\,--\,(#5#2#6)}
\crefmultiformat{equation}{(#2#1#3)}{\,--\,(#2#1#3)}{#2#1#3}{#2#1#3}

\crefname{equation}{}{}
\crefname{figure}{Figure}{Figures}
\crefname{algorithm}{Algorithm}{Algorithms}
\crefname{table}{Table}{Tables}
\crefname{lemma}{Lemma}{Lemmas}
\crefname{theorem}{Theorem}{Theorems}
\crefname{section}{Section}{Sections}
\crefname{definition}{Definition}{Definitions}
\commentfalse
\IEEEoverridecommandlockouts                              
\title{\LARGE \bf Convex Relaxation of Bilinear Matrix Inequalities\\Part I: Theoretical Results}
\author{Mohsen Kheirandishfard, Fariba Zohrizadeh, and Ramtin Madani
\thanks{Mohsen Kheirandishfard and Fariba Zohrizadeh are with the Department of Computer Science and Engineering, The University of Texas at Arlington, Arlington, TX 76019, USA (emails:mohsen.kheirandishfard@uta.edu, fariba.zohrizadeh@uta.edu), Ramtin Madani is with the Department of Electrical Engineering, The University of Texas at Arlington, Arlington, TX 76019, USA (email:ramtin.madani@uta.edu). This work is in part supported by the NSF award 1809454 and a University of Texas System STARs award.}%
}

\begin{document}

\maketitle
\thispagestyle{empty}
\pagestyle{empty}

\begin{abstract}
This two-part paper is concerned with the problem of minimizing a linear objective function subject to a bilinear matrix inequality (BMI) constraint. In this part, we first consider a family of convex relaxations which transform BMI optimization problems into polynomial-time solvable surrogates. As an alternative to the state-of-the-art semidefinite programming (SDP) and second-order cone programming (SOCP) relaxations, a computationally efficient {\it parabolic relaxation} is developed, which relies on convex quadratic constraints only. Next, we developed a family of penalty functions, which can be incorporated into the objective of SDP, SOCP, and parabolic relaxations to facilitate the recovery of feasible points for the original non-convex BMI optimization. Penalty terms can be constructed using any arbitrary initial point. We prove that if the initial point is sufficiently close to the feasible set, then the penalized relaxations are guaranteed to produce feasible points for the original BMI. In Part II of the paper, the efficacy of the proposed penalized convex relaxations is demonstrated on benchmark instances of $\mathcal{H}_2$ and $\mathcal{H}_{\infty}$ optimal control synthesis problems.
\end{abstract}
\section{INTRODUCTION}
A wide range of real-world problems in the area of control can be cast as optimization problems with matrix inequality constraints \cite{el2000advances,boyd1994linear,gahinet1994linear,kazemi2013adaptive}. As a special case, the class of problems with linear matrix inequalities (LMIs) can be solved efficiently up to any desired accuracy via interior-point based methods \cite{gahinet1993general,nesterov1994interior,ben2001lectures,vandenberghe1997algorithms}. However, despite various applications, optimization in the presence of bilinear matrix inequalities (BMIs) is computationally prohibitive and NP-hard in general. Significant efforts have been devoted to the development of algorithms 
for solving BMI optimization problems \cite{kovcvara2005nonlinear,ma2012multi,correa2004global}, including software packages \cite{orsi2006newton,fiala2013penlab,kocvara2005penbmi}.


In \cite{goh1996robust,safonov1994control,safonov1994diameter,goh1994global}, a number of robust control design problems have been formulated using BMI constraints. Later, the use of BMIs has been extended to other control applications, such as state- and output-feedback controller design \cite{d2003distributed,fattahi2016transformation,wang2016feasibility,nian2013bilinear,hassibi1998low,chamanbaz2014statistical,fazelnia2017convex}, affine fuzzy system design \cite{feng2005h,qiu2013static}, stability of fractional-order systems \cite{lim2013stability}. Despite valuable insights gained from BMI formulations of control-theoretic problems, the design of an efficient algorithm for solving the resulting non-convex formulations has remained an open problem \cite{vanantwerp2000tutorial}. In \cite{toker1995np,blondel1997np}, it is proven that a general BMI constrained optimization problem is NP-hard.

Despite the above-mentioned theoretical barrier, various approaches have been developed in the literature to tackle BMI optimization problems. In \cite{goh1994biaffine,el1994synthesis}, alternating minimization (AM)-based algorithms are proposed which divide variables into two blocks that can be alternately optimized until convergence. Although AM-based methods enjoy simple implementation and perform satisfactorily in many cases, they offer no convergence guarantees to a feasible solution. 
Another approach is to solve a sequence of convex relaxations until a satisfactory solution is obtained
\cite{liu1996numerical,hassibi1999path,safonov1994diameter,ibaraki2001rank,doelman2016sequential,beck2010sequential}. In \cite{lee2016sequential,hassibi1999path}, BMI optimization problems are tackled by forming a sequence of semidefinite programming (SDP) relaxations. In \cite{dinh2012combining}, a sequential method is developed based on difference-of-convex programming with convergence guarantees to (sub)-optimal solutions. In \cite{ibaraki2001rank,doelman2016sequential,el1997cone} rank-constrained formulations with nuclear norm penalties are investigated along with bound-tightening methods for solving general BMI optimization problems. In \cite{fukuda2001branch,goh1994global,tuan1999new,tuan2000low}, branch-and-bound (BB) methods are developed 
with convergence guarantees to global optimality.
In general, BB methods are computationally prohibitive and their applicability is limited to moderate-sized problems. A novel global optimization method has been recently presented in \cite{chiu2017method} which 
tackles BMI problems using hybrid multi-objective optimization methods. 
BMIs can be categorized as a special case of polynomial matrix inequalities. Therefore, methods for solving general polynomial matrix inequalities are applicable to BMIs as well \cite{kim2005generalized,ahmadi2015sum}. Despite computational complexity for real-world applications, the most notable example is Lasserre's hierarchy of LMI relaxations \cite{lasserre2001global}, based on which several software packages have been developed \cite{prajna2004sostools,papachristodoulou2013sostools,henrion2003gloptipoly}.

The main contribution of this work is a novel and general convex relaxation, regarded as parabolic relaxation, for solving optimization problems with BMI constraints. The proposed convex relaxation relies on convex quadratic constraints as opposed to the SDP and SOCP relaxations which rely on computationally expensive conic constraints.
Our second contribution is concerned with finding feasible and near-globally optimal solutions for BMI optimization problems. To this end, we incorporate a penalty term into the objective function of convex relaxations. The proposed penalty term is compatible with SDP, SOCP, and parabolic relaxations, and can be customized using any available initial point. We prove that If the initial point is feasible for the original BMI problem, then the outcome of penalized relaxation is guaranteed to be feasible as well. Moreover, 
any infeasible initial point which is close to the BMI feasible set is guaranteed to produce feasible point. 
In Part II, the proposed penalized relaxation framework is extended to a sequential scheme, which is tested on a variety of optimal control design problems.

\subsection{Notation}
Throughout the paper, the scalars, vectors, and matrices are respectively shown by italic letters, lower-case bold letters, and upper-case bold letters. Symbols $\Rbb$, $\Rbb^{n}$, and $\Rbb^{n\times m}$ respectively denote the set of real scalars, real vectors of size $n$, and real matrices of size $n\times m$. The set of real $n\times n$ symmetric matrices and positive semidefinite matrices are shown with $\Sbb_n$ and $\Sbb_n^+$, respectively. For given vector $\abf$ and matrix $\Abf$, symbols $a_i$ and $A_{ij}$ respectively indicate the $i^{th}$ element of $\abf$ and $(i,j)^{th}$ element of $\Abf$. Notations $[a]_{i\in\mathcal{I}}$ and $[A]_{ij\in\mathcal{I}}$ respectively shows the sub-vector and sub-matrix corresponding to the set of indices $\mathcal{I}$. Notation $\Abf\succeq 0$ means $\Abf$ is positive-semidefinite ($\Abf\succ 0$ indicates positive definite) and $\Abf\preceq 0$ means $\Abf$ is negative-semidefinite ($\Abf\prec 0$ indicates negative definite). For two given matrices $\Abf$ and $\Bbf$ of the same size, symbol $\langle\Abf,\Bbf\rangle=\tr\{\Abf^{\!\top}\Bbf\}$ shows the inner product between the matrices where $\tr\{.\}$ and $(.)^\top$ respectively denote the trace and transpose operators. Notation $\lVert.{\rVert}_p$ refers to either matrix norm or vector norm depending on the context and $\lvert.\rvert$ indicates the absolute value. Symbols $\Ibf$, $\ebf_i$, and $\mathbf{0}$ denote the identity matrix, standard basis vector, and zero matrix of appropriate dimensions, respectively. Letters $\Ncal$ and $\Mcal$ are shorthand for sets $\{1,\dots,n\}$ and $\{1,\dots,m\}$, respectively.

\section{Problem Formulation}
This paper is concerned with the following class of optimization problems with linear objective and a bilinear matrix inequality (BMI) constraint:
\begin{subequations} \label{eq:BMI_BMIOrig}
\begin{align}
&\underset{{\xbf}\in{\Rbb}^{n}}{\text{minimize}} && \cbf^{\!\top}\xbf \label{eq:BMI_BMIOrig_obj}\\
&\text{subject to}
&&p(\xbf,\xbf\xbf^{\top})\preceq 0,\label{eq:BMI_BMIOrig_con_01}
\end{align}
\end{subequations}
where $\cbf\in\Rbb^{n}$ is the cost vector, and $p\!:\!\Rbb^{n}\times \Sbb_{n}\!\rightarrow\Sbb_{m}$ is a linear matrix-valued function, which is regarded as matrix pencil. In general, $p$ can be formulated as:
\begin{equation}\label{eq:BMI_MatPen}
\begin{aligned}[b]
& p(\xbf,\Xbf)\triangleq{\Fbf}_{\sm{0}{6}}+\sum_{\sm{k\in\Ncal}{5.5}} {x}_{k}{\Kbf}_k+\sum_{\sm{i\in\Ncal}{5.5}}\sum_{\sm{j\in\Ncal}{5.5}}{X}_{ij}{\Lbf}_{ij}.
\end{aligned}
\end{equation}
where ${\Fbf_0}$, $\{{\Kbf}_k\}_{\sm{k\in\Ncal}{5.5}}$, and $\{{\Lbf}_{ij}\}_{\sm{i,j\in\Ncal}{5.5}}$ are $m\times m$ real symmetric matrices. With no loss of generality, we can assume that $\Lbf_{ij}=\Lbf_{ji}$ for all $i,j\!\in\!\Ncal$, since $\Xbf$ is a symmetric matrix.

Problem \cref{eq:BMI_BMIOrig_obj,eq:BMI_BMIOrig_con_01} is non-convex and NP-hard in general, due to the presence of the BMI constraint \cref{eq:BMI_BMIOrig_con_01}. To tackle this problem, it is common practice to solve convex surrogates that produce lower bounds on the globally-optimal cost of the original non-convex problem \cref{eq:BMI_BMIOrig_obj,eq:BMI_BMIOrig_con_01}. To this end, an auxiliary matrix variable $\Xbf$ is introduced to account for ${\xbf}{\xbf}^{\top}$. This leads to the following {\it lifted reformulation} of the problem \cref{eq:BMI_BMIOrig_obj,eq:BMI_BMIOrig_con_01}:
\begin{subequations} \label{eq:BMI_BMIOrig_lifted}
\begin{align}
&\underset{\begin{subarray}{l} {\xbf}\in{\Rbb}^{n}\!,{\Xbf}\in{\Sbb}_{n}\end{subarray}}{\text{minimize}} && \cbf^{\!\top}\xbf\label{eq:BMI_BMIOrig_lifted_obj}\\
&\text{subject to}
&&p(\xbf,\Xbf)\preceq 0,\label{eq:BMI_BMIOrig_lifted_con_01} \\
&&& {\Xbf}={\xbf}{\xbf}^{\top}, \label{eq:BMI_BMIOrig_lifted_con_02}
\end{align}
\end{subequations}
where constraint \cref{eq:BMI_BMIOrig_lifted_con_02} is imposed to preserve the equivalency. 
Lifting casts the problem into a higher dimensional space in which the BMI constraint \eqref{eq:BMI_BMIOrig_con_01} is transformed into a linear matrix inequality and the entire non-convexity is captured by the new constraint \cref{eq:BMI_BMIOrig_lifted_con_02}. In what follows, we will substitute \cref{eq:BMI_BMIOrig_lifted_con_02} with convex alternatives and revise the objective function in order to obtain feasible and near-globally optimal points for the original problem \cref{eq:BMI_BMIOrig_obj,eq:BMI_BMIOrig_con_01}.

\section{Preliminaries}
In order to further analyze the BMI constraint \eqref{eq:BMI_BMIOrig_con_01}, the next definition introduces the notion of pencil norm.

\vspace{1.5mm}
\begin{definition}[Pencil Norm]
	For every $q\geq1$, the $q$-norm of the matrix pencil in equation \eqref{eq:BMI_MatPen} is defined as
	\begin{equation}\label{eq:MOR_DefPencNorm}
	\begin{aligned}[b] 
	{\lVert p\rVert}_{q}\!\triangleq\msh\text{\normalfont max}\Big\{\big\lVert\msh\big[\ubf^{\!\top}\msh\Lbf_{ij}\ubf{\big]}_{\sm{i,j\in\Ncal^2}{5.5}}{\big\rVert}_{q}\Big|\forall \ubf\msh\in\msh\Rbb^m,\lVert\ubf{\rVert}_{2}\!=\!1\Big\}.
	\end{aligned}
	\end{equation}
\end{definition}
\vspace{1.5mm}

The next definition provides a measure of distance between any arbitrary point in $\Rbb^n$ and the feasible set of optimization problem \cref{eq:BMI_BMIOrig_obj,eq:BMI_BMIOrig_con_01}. 
\vspace{1.5mm}
\begin{definition}[Feasiblity Distance]\label{def:DistFromSet}
	For every $\xbf\in\Rbb^n$, define the feasibility distance $d_{\Fcal}:\Rbb^n\rightarrow\Rbb$ as 
	\begin{equation}\label{eq:BMI_DefDistFromSet}
	\begin{aligned}
	d_{\sm{\Fcal}{5.5}}(\xbf)\triangleq\text{\normalfont inf}\{\lVert\xbf-\abf{\rVert}_{2}\;\big|\;\abf\in\Fcal\},
	\end{aligned}
	\end{equation}
	where $\Fcal\subseteq\Rbb^{n}$ denotes the feasible set of the BMI problem \cref{eq:BMI_BMIOrig_obj,eq:BMI_BMIOrig_con_01}. Observe that the feasibility distance is equal to 0 if $\xbf\in\Fcal$.
\end{definition}
\vspace{1.5mm}
We use the Mangasarian-Fromovitz constraint qualification (MFCQ) condition from \cite{shapiro1997first} in order to characterize well-behaved feasible points of problem \cref{eq:BMI_BMIOrig_obj,eq:BMI_BMIOrig_con_01}. 
\vspace{1.5mm}
\begin{definition}[MFCQ Condition]\label{def:BMI_MFCQCond}
	A feasible point $\xbf\in\mathcal{F}$ of problem \cref{eq:BMI_BMIOrig_obj,eq:BMI_BMIOrig_con_01} is said to
	satisfy the MFCQ condition if there exists $\bbf\in\Rbb^{n}$ such that
	\begin{equation}\label{eq:BMI_AssumpMFCond}
 p(\xbf,\xbf{\xbf}^{\!\top})+\sum_{\sm{k\in\Ncal}{5.5}} b_{k}(\Kbf_{k}\!+\!{\delta}_{k}(\xbf))\prec 0,
	\end{equation}
	where for every $k\in\Ncal$, the matrix function $\delta_{k}:\Rbb^{n}\rightarrow\Sbb_m$ is defined as
	\begin{equation} \label{eq:BMI_DefDelta}
	\delta_{k}(\xbf)\triangleq 2\sum_{\sm{i\in\Ncal}{5.5}} x_{i}\Lbf_{ki},
	\end{equation}
	representing the derivative of bilinear terms of pencil $p$ with respect to $x_k$. 
\end{definition}
\vspace{1.5mm}

In the following definition, we introduce a generalization of the MFCQ condition to cover infeasible points as well.
\vspace{1.5mm}
\begin{definition}[G-MFCQ Condition]
	An arbitrary point $\xbf\in\mathbb{R}^n$ is said to satisfy the Generalized Mangasarian-Fromovitz constraint qualification (G-MFCQ) condition for problem \cref{eq:BMI_BMIOrig_obj,eq:BMI_BMIOrig_con_01}, if there exists $\bbf\in\Rbb^{n}$ where
	\begin{equation}\label{eq:BMI_AssumpGMFCond}
	\begin{aligned}[b]
	\sum_{\sm{k\in\Ncal}{5.5}} b_{k}(\Kbf_{k}\!+\!{\delta}_{k}(\xbf))\prec 0.
	\end{aligned}
	\end{equation}
Moreover, define the G-MFCQ function $s:\Rbb^{n}\rightarrow\Rbb$ as
	\begin{equation} 
	\begin{aligned}
	&s(\xbf)\triangleq{\text{\normalfont max}}\Big\{\ubar{\lambda}\Big(\!-\!\sum_{\sm{k\in\Ncal}{5.5}}b_{k}(\Kbf_{k}\!+\!{\delta}_{k}(\xbf)\big)\!\Big)\Big|\lVert\bbf{\rVert}_{2}\!=\!1\Big\}, \label{eq:MOR_DefsXGMF}
	\end{aligned}
	\end{equation}
	where the operator $\ubar{\lambda}(.)$ returns the minimum eigenvalue of its input argument.
\end{definition}
\vspace{1.5mm}
\section{Convex Relaxation}
This section aims at introducing a family of convex relaxations for the lifted problem \cref{eq:BMI_BMIOrig_lifted_obj,eq:BMI_BMIOrig_lifted_con_01,eq:BMI_BMIOrig_lifted_con_02}. 
Consider the following formulation: 
\begin{subequations} \label{eq:BMI_GenRelaxation}
\begin{align}
&\underset{{\xbf}\sm{\in}{6}{\Rbb}^{n}\!,{\Xbf}\sm{\in}{6}{\Sbb}_{n}}{\text{minimize}} && \cbf^{\!\top} \xbf \label{eq:BMI_GenRelaxation_obj}\\
&\text{subject to}
&&p(\xbf,\Xbf)\;\preceq 0,\label{eq:BMI_GenRelaxation_con_01} \\
&&& \Xbf-\xbf\xbf^{\!\top}\!\in\Ccal, \label{eq:BMI_GenRelaxation_con_02}
\end{align}
\end{subequations}
in which $\Ccal\subseteq\mathbb{S}_n$. Observe that the problems \cref{eq:BMI_GenRelaxation_obj,eq:BMI_GenRelaxation_con_01,eq:BMI_GenRelaxation_con_02} and \cref{eq:BMI_BMIOrig_lifted_obj,eq:BMI_BMIOrig_lifted_con_01,eq:BMI_BMIOrig_lifted_con_02} are equivalent if $\Ccal=\{0\}$. 

We consider different choices for $\Ccal$, which make the constraint \eqref{eq:BMI_GenRelaxation_con_02} convex. First, the standard semidefinite programming (SDP) and second-order cone programming (SOCP) relaxations are discussed and then, we introduce a novel {\it parabolic relaxation}, which transforms the constraint \cref{eq:BMI_BMIOrig_lifted_con_02} into a set of convex quadratic inequalities. The optimal cost for each of the above convex relaxation can serve as a lower bound for the global cost of the original problem \cref{eq:BMI_BMIOrig_obj,eq:BMI_BMIOrig_con_01}. If the optimal solution of a relaxed problem satisfies \cref{eq:BMI_BMIOrig_lifted_con_02}, the relaxation is regarded as {\it exact}. 


\subsection{Semidefinite Programming Relaxation}
The following choice for $\Ccal$ leads to the SDP relaxation of the problem \cref{eq:BMI_GenRelaxation_obj,eq:BMI_GenRelaxation_con_01,eq:BMI_GenRelaxation_con_02}:
\begin{equation}\label{eq:BMI_SDP_Set}
\Ccal_{\sm{\mathrm{1}}{5}}\!=\!\{\Hbf\in\Sbb_{n}\!\mid\!\Hbf\succeq 0\}.
\end{equation}
If $\Ccal=\Ccal_{\sm{\mathrm{1}}{5}}$ the optimization problem \cref{eq:BMI_GenRelaxation_obj,eq:BMI_GenRelaxation_con_01,eq:BMI_GenRelaxation_con_02} boils down to a semidefinite program, which can be efficiently solved in polynomial time up to any desired accuracy using the existing methods. 
\subsection{Second-Order Cone Programming Relaxation}
Semidefinite programming can be computationally demanding and its application is limited to small-scale problems. A popular alternative is the SOCP relaxation 
which can be deduced from the following choice for $\Ccal$:
\begin{equation}\label{eq:BMI_SOCP_Set}
\begin{aligned}[b]
\!\!\!\!\Ccal_{\sm{\mathrm{2}}{5}}\!=\!\{\Hbf\in\Sbb_{n} \mid H_{ii}\geq\! 0,\;\;H_{ii}H_{jj}\geq H_{ij}^2, \;\;\sm{\forall i,\!j\!\in\Ncal}{8}\}.
\end{aligned}
\end{equation}
It is straightforward to show that $\Ccal_{\sm{\mathrm{1}}{5}}$ is a subset of $\Ccal_{\sm{\mathrm{2}}{5}}$, which implies that the lower bounds from SDP relaxation are guaranteed to be tighter than or equal to the lower bounds obtained by SOCP relaxation.
\subsection{Parabolic Relaxation}
In this subsection, the parabolic relaxation is introduced as a computationally efficient alternative to SDP and SOCP relaxations. Parabolic relaxation transforms the non-convex constraint \cref{eq:BMI_BMIOrig_lifted_con_02} to a number of convex quadratic inequalities. To formulate the parabolic relaxation of the problem \cref{eq:BMI_GenRelaxation_obj,eq:BMI_GenRelaxation_con_01,eq:BMI_GenRelaxation_con_02}, the following choice for $\Ccal$ should be employed:
\begin{equation}\label{eq:BMI_QP_Set}
\begin{aligned}[b]
\!\!\!\!\Ccal_{\sm{\mathrm{3}}{5}}\!=\!\{\Hbf\in\Sbb_{n}\!\mid\! H_{ii}\!\geq\! 0,H_{ii}\!+\!H_{jj}\!\geq\! 2\left|H_{ij}\right|\!,\sm{\forall i,\!j\!\in\Ncal}{8}\}.
\end{aligned}
\end{equation}
It can be easily observed that if $\Ccal=\Ccal_{\sm{\mathrm{3}}{5}}$, then the constraint \eqref{eq:BMI_GenRelaxation_con_02} is equivalent to the following quadratic inequalities:
\begin{subequations}
\begin{align}
X_{ii}+X_{jj}-2X_{ij} \geq (x_i-x_j)^2\qquad \forall i,\!j\!\in\Ncal,   \\
X_{ii}+X_{jj}+2X_{ij} \geq (x_i+x_j)^2\qquad \forall i,\!j\!\in\Ncal,
\end{align}
\end{subequations}
which means that the parabolic relaxation is computationally cheaper than the SDP and SOCP relaxations.

Note that the presented relaxations are not necessarily exact. In the next section, the objective function \cref{eq:BMI_GenRelaxation_obj} is revised to facilitate the recovery of feasible points for the original non-convex problem \cref{eq:BMI_BMIOrig_obj,eq:BMI_BMIOrig_con_01}.

\section{Penalized Convex Relaxation} \label{sec:pana}
The {\it penalized convex relaxation} of the BMI optimization \cref{eq:BMI_BMIOrig_obj,eq:BMI_BMIOrig_con_01} is given as
\begin{subequations} \label{eq:BMI_GenRelaxation_Pen}
\begin{align}
&\underset{{\xbf}\sm{\in}{6}{\Rbb}^{n}\!,{\Xbf}\sm{\in}{6}{\Sbb}_{n}}{\text{minimize}} && \cbf^{\!\top} \xbf+{\eta}\;(\mathrm{tr}\{\Xbf\}-2\;{\check{\xbf}}^{\!\top}{\xbf}+\check{\xbf}^{\!\top}{\check{\xbf}}) \label{eq:BMI_GenRelaxation_Pen_obj}\\
&\text{subject to}
&&p(\xbf,\Xbf)\;\preceq 0,\label{eq:BMI_GenRelaxation_Pen_con_01} \\
&&& \Xbf-\xbf\xbf^{\!\top}\!\in\Ccal, \label{eq:BMI_GenRelaxation_Pen_con_02}
\end{align}
\end{subequations}
where $\check{\xbf}\in\mathbb{R}^n$ is an initial guess for the unknown solution (either feasible or infeasible), $\eta>0$ is a regularization parameter, which offers a trade-off between the original objective function and the penalty term, and $\Ccal\in\{\Ccal_{\sm{\mathrm{1}}{5}},\Ccal_{\sm{\mathrm{2}}{5}},\Ccal_{\sm{\mathrm{3}}{5}}\}$.
\vspace{1.5mm}

The next theorem states that if the initial point $\check{\xbf}$ is feasible and satisfies MFCQ, then the penalized convex relaxation preserves the feasibility of $\check{\xbf}$ and produces a solution with improved objective value.

\vspace{1.5mm}
\begin{theorem}\label{thm:BMI_thm_01}
Assume that $\check{\xbf}\in\mathcal{F}$ is a feasible point for problem \cref{eq:BMI_BMIOrig_obj,eq:BMI_BMIOrig_con_01} that satisfies the MFCQ condition. If $\Ccal\in\{\Ccal_{\sm{\mathrm{1}}{5}},\Ccal_{\sm{\mathrm{2}}{5}},\Ccal_{\sm{\mathrm{3}}{5}}\}$ and $\eta$ is sufficiently large, then the penalized convex relaxation problem \cref{eq:BMI_GenRelaxation_Pen_obj,eq:BMI_GenRelaxation_Pen_con_02} has a unique solution $(\accentset{\ast}{\xbf},\accentset{\ast}{\Xbf})$, which satisfies $\accentset{\ast}{\Xbf}=\accentset{\ast}{\xbf}\accentset{\ast}{\xbf}^{\top}$ and $\cbf^{\top}\accentset{\ast}{\xbf}\leq\cbf^{\top}\check{\xbf}$.
\end{theorem}
\vspace{1.5mm}
\begin{proof}
See Appendix for the proof.
\end{proof}
\vspace{1.5mm}
According to \cref{thm:BMI_thm_01}, the proposed penalized relaxation preserves the feasibility of the initial point. In what follows, we show that if the initial point is not feasible for \cref{eq:BMI_BMIOrig_obj,eq:BMI_BMIOrig_con_01}, but sufficiently close to its feasible set, then the penalized convex relaxation problem \cref{eq:BMI_GenRelaxation_Pen_obj,eq:BMI_GenRelaxation_Pen_con_02} is guaranteed to produce a feasible solution as well.
\vspace{1.5mm}
\begin{theorem}\label{thm:BMI_thm_02}
Assume that $k\in\{1,2,3\}$ and $\Ccal=\Ccal_k$.
Consider an arbitrary point $\check{\xbf}\in\Rbb^n$, which satisfies the G-MFCQ condition for problem \cref{eq:BMI_BMIOrig_obj,eq:BMI_BMIOrig_con_01}, and let 
\begin{align}
\frac{d_{\sm{\Fcal}{5.5}}(\check{\xbf})}{s(\check{\xbf})} \leq \frac{\omega_k}{\lVert p{\rVert}_2}
\end{align}
where $\omega_1=4^{-1}$, $\omega_2=(2n)^{-1}$, and $\omega_3=(2+2\sqrt{n})^{-1}$.
If $\eta$ is sufficiently large, then the penalized convex relaxation problem \cref{eq:BMI_GenRelaxation_Pen_obj,eq:BMI_GenRelaxation_Pen_con_02} 
has a unique solution $(\accentset{\ast}{\xbf},\accentset{\ast}{\Xbf})$, which satisfies $\accentset{\ast}{\Xbf}=\accentset{\ast}{\xbf}\accentset{\ast}{\xbf}^{\top}$.
\end{theorem}
\vspace{1.5mm}
\begin{proof}
See Appendix for the proof.
\end{proof}
\vspace{1.5mm}

In Part II of this paper, we use the results of Theorems 1 and 2 to developed a sequential method for solving general BMI optimization problems. 

\section{Conclusions}
In this paper, a variety of convex relaxation methods are introduced for solving the class of optimization problems with bilinear matrix inequality (BMI) constraints. First, the well-known SDP and SOCP relaxations are discussed, and then a novel parabolic relaxation is introduced as a low-complexity alternative to conic relaxations. We propose a penalization method which is compatible with SDP, SOCP, and parabolic relaxations, and is able to produce feasible solutions for the original non-convex BMI optimization problem. 
In part II, the proposed penalized convex relaxation scheme is generalized to a sequential framework and is tested on a variety of challenging benchmark optimal control design problems.

\appendices
\bibliographystyle{IEEEtran}
\bibliography{IEEEabrv,egbib}

\begin{thebibliography}{10}
\providecommand{\url}[1]{#1}
\csname url@samestyle\endcsname
\providecommand{\newblock}{\relax}
\providecommand{\bibinfo}[2]{#2}
\providecommand{\BIBentrySTDinterwordspacing}{\spaceskip=0pt\relax}
\providecommand{\BIBentryALTinterwordstretchfactor}{4}
\providecommand{\BIBentryALTinterwordspacing}{\spaceskip=\fontdimen2\font plus
\BIBentryALTinterwordstretchfactor\fontdimen3\font minus
  \fontdimen4\font\relax}
\providecommand{\BIBforeignlanguage}[2]{{%
\expandafter\ifx\csname l@#1\endcsname\relax
\typeout{** WARNING: IEEEtran.bst: No hyphenation pattern has been}%
\typeout{** loaded for the language `#1'. Using the pattern for}%
\typeout{** the default language instead.}%
\else
\language=\csname l@#1\endcsname
\fi
#2}}
\providecommand{\BIBdecl}{\relax}
\BIBdecl

\bibitem{el2000advances}
L.~El~Ghaoui and S.-l. Niculescu, \emph{Advances in linear matrix inequality
  methods in control}.\hskip 1em plus 0.5em minus 0.4em\relax SIAM, 2000.

\bibitem{boyd1994linear}
S.~Boyd, L.~El~Ghaoui, E.~Feron, and V.~Balakrishnan, \emph{Linear matrix
  inequalities in system and control theory}.\hskip 1em plus 0.5em minus
  0.4em\relax SIAM, 1994.

\bibitem{gahinet1994linear}
P.~Gahinet and P.~Apkarian, ``A linear matrix inequality approach to
  {$\Hcal_{\infty}$} control,'' \emph{INT J ROBUST NONLIN}, vol.~4, no.~4, pp.
  421--448, 1994.

\bibitem{kazemi2013adaptive}
H.~Kazemi and M.~Namvar, ``Adaptive compensation of actuator dynamics in
  manipulators without joint torque measurement,'' in \emph{Proc. IEEE Conf.
  Decis. Control}.\hskip 1em plus 0.5em minus 0.4em\relax IEEE, 2013, pp.
  2294--2299.

\bibitem{gahinet1993general}
P.~Gahinet and A.~Nemirovskii, ``General-purpose {LMI} solvers with
  benchmarks,'' in \emph{Proc. IEEE Conf. Decis. Control}, 1993, pp.
  3162--3165.

\bibitem{nesterov1994interior}
Y.~Nesterov and A.~Nemirovskii, \emph{Interior-point polynomial algorithms in
  convex programming}.\hskip 1em plus 0.5em minus 0.4em\relax SIAM, 1994.

\bibitem{ben2001lectures}
A.~Ben-Tal and A.~Nemirovski, \emph{Lectures on modern convex optimization:
  analysis, algorithms, and engineering applications}.\hskip 1em plus 0.5em
  minus 0.4em\relax SIAM, 2001.

\bibitem{vandenberghe1997algorithms}
L.~Vandenberghe and V.~Balakrishnan, ``Algorithms and software for {LMI}
  problems in control,'' \emph{IEEE. Control. Syst.}, vol.~17, no.~5, pp.
  89--95, 1997.

\bibitem{kovcvara2005nonlinear}
M.~Ko{\v{c}}vara, F.~Leibfritz, M.~Stingl, and D.~Henrion, ``A nonlinear {SDP}
  algorithm for static output feedback problems in {COMPl$_\mathrm{e}$ib},''
  \emph{IFAC-PapersOnLine}, vol.~38, no.~1, pp. 1055--1060, 2005.

\bibitem{ma2012multi}
L.~Ma, X.~Meng, Z.~Liu, and L.~Du, ``Multi-objective and reliable control for
  trajectory-tracking of rendezvous via parameter-dependent {L}yapunov
  functions,'' \emph{Acta Astronaut.}, vol.~81, no.~1, pp. 122--136, 2012.

\bibitem{correa2004global}
R.~Correa, ``A global algorithm for nonlinear semidefinite programming,''
  \emph{SIAM J. Optim.}, vol.~15, no.~1, pp. 303--318, 2004.

\bibitem{orsi2006newton}
R.~Orsi, U.~Helmke, and J.~B. Moore, ``A {N}ewton-like method for solving rank
  constrained linear matrix inequalities,'' \emph{Automatica}, vol.~42, no.~11,
  pp. 1875--1882, 2006.

\bibitem{fiala2013penlab}
J.~Fiala, M.~Ko{\v{c}}vara, and M.~Stingl, ``{PENLAB}: A {MATLAB} solver for
  nonlinear semidefinite optimization,'' \emph{arXiv preprint arXiv:1311.5240},
  2013.

\bibitem{kocvara2005penbmi}
M.~Kocvara, M.~Stingl, and P.~GbR, ``{PENBMI} user’s guide (version 2.0),''
  \emph{software manual, PENOPT GbR, Hauptstrasse A}, vol.~31, p. 91338, 2005.

\bibitem{goh1996robust}
K.-C. Goh, M.~Safonov, and J.~Ly, ``Robust synthesis via bilinear matrix
  inequalities,'' \emph{INT J ROBUST NONLIN}, vol.~6, no. 9-10, pp. 1079--1095,
  1996.

\bibitem{safonov1994control}
M.~G. Safonov, K.-C. Goh, and J.~Ly, ``Control system synthesis via bilinear
  matrix inequalities,'' in \emph{Proc. IEEE Am. Control Conf.}, vol.~1, 1994,
  pp. 45--49.

\bibitem{safonov1994diameter}
M.~Safonov and G.~Papavassilopoulos, ``The diameter of an intersection of
  ellipsoids and {BMI} robust synthesis,'' in \emph{Proc. of the IFAC Symp. on
  Robust Control Design}, 1994, pp. 313--317.

\bibitem{goh1994global}
K.-C. Goh, M.~Safonov, and G.~Papavassilopoulos, ``A global optimization
  approach for the {BMI} problem,'' in \emph{Proc. IEEE Conf. Decis. Control},
  vol.~3, 1994, pp. 2009--2014.

\bibitem{d2003distributed}
R.~D'Andrea and G.~E. Dullerud, ``Distributed control design for spatially
  interconnected systems,'' \emph{IEEE Trans. Autom. Control}, vol.~48, no.~9,
  pp. 1478--1495, 2003.

\bibitem{fattahi2016transformation}
S.~Fattahi, G.~Fazelnia, and J.~Lavaei, ``Transformation of optimal centralized
  controllers into near-globally optimal static distributed controllers,''
  \emph{IEEE Trans. Autom. Control}, 2017.

\bibitem{wang2016feasibility}
Y.~Wang and R.~Rajamani, ``Feasibility analysis of the bilinear matrix
  inequalities with an application to multi-objective nonlinear observer
  design,'' in \emph{Proc. IEEE Conf. Decis. Control}, 2016, pp. 3252--3257.

\bibitem{nian2013bilinear}
X.~Nian, Z.~Sun, H.~Wang, H.~Zhang, and X.~Wang, ``Bilinear matrix inequality
  approaches to robust guaranteed cost control for uncertain discrete-time
  delay system,'' \emph{OPTIM CONTR APPL MET}, vol.~34, no.~4, pp. 433--441,
  2013.

\bibitem{hassibi1998low}
A.~Hassibi, J.~How, and S.~Boyd, ``Low-authority controller design via convex
  optimization,'' in \emph{Proc. IEEE Conf. Decis. Control}, vol.~1, 1998, pp.
  140--145.

\bibitem{chamanbaz2014statistical}
M.~Chamanbaz, F.~Dabbene, R.~Tempo, V.~Venkataramanan, and Q.-G. Wang, ``A
  statistical learning theory approach for uncertain linear and bilinear matrix
  inequalities,'' \emph{Automatica}, vol.~50, no.~6, pp. 1617--1625, 2014.

\bibitem{fazelnia2017convex}
G.~Fazelnia, R.~Madani, A.~Kalbat, and J.~Lavaei, ``Convex relaxation for
  optimal distributed control problems,'' \emph{IEEE Trans. Autom. Control},
  vol.~62, no.~1, pp. 206--221, 2017.

\bibitem{feng2005h}
G.~Feng, C.-L. Chen, D.~Sun, and Y.~Zhu, ``{$\Hcal_{\infty}$} controller
  synthesis of fuzzy dynamic systems based on piecewise {L}yapunov functions
  and bilinear matrix inequalities,'' \emph{IEEE Trans. Fuzzy Syst}, vol.~13,
  no.~1, pp. 94--103, 2005.

\bibitem{qiu2013static}
J.~Qiu, G.~Feng, and H.~Gao, ``Static-output-feedback {$\Hcal_{\infty}$}
  control of continuous-time {T}-{S} fuzzy affine systems via piecewise
  {L}yapunov functions,'' \emph{IEEE Trans. Fuzzy Syst}, vol.~21, no.~2, pp.
  245--261, 2013.

\bibitem{lim2013stability}
Y.-H. Lim, K.-K. Oh, and H.-S. Ahn, ``Stability and stabilization of
  fractional-order linear systems subject to input saturation,'' \emph{IEEE
  Trans. Autom. Control}, vol.~58, no.~4, pp. 1062--1067, 2013.

\bibitem{vanantwerp2000tutorial}
J.~G. VanAntwerp and R.~D. Braatz, ``A tutorial on linear and bilinear matrix
  inequalities,'' \emph{J Process Control}, vol.~10, no.~4, pp. 363--385, 2000.

\bibitem{toker1995np}
O.~Toker and H.~Ozbay, ``On the {NP}-hardness of solving bilinear matrix
  inequalities and simultaneous stabilization with static output feedback,'' in
  \emph{Proc. IEEE Am. Control Conf.}, vol.~4, 1995, pp. 2525--2526.

\bibitem{blondel1997np}
V.~Blondel and J.~N. Tsitsiklis, ``{NP}-hardness of some linear control design
  problems,'' \emph{SIAM J. Control Optim.}, vol.~35, no.~6, pp. 2118--2127,
  1997.

\bibitem{goh1994biaffine}
K.~Goh, L.~Turan, M.~Safonov, G.~Papavassilopoulos, and J.~Ly, ``Biaffine
  matrix inequality properties and computational methods,'' in \emph{Proc. IEEE
  Am. Control Conf.}, vol.~1, 1994, pp. 850--855.

\bibitem{el1994synthesis}
L.~El~Ghaoui and V.~Balakrishnan, ``Synthesis of fixed-structure controllers
  via numerical optimization,'' in \emph{Proc. IEEE Conf. Decis. Control},
  vol.~3, 1994, pp. 2678--2683.

\bibitem{liu1996numerical}
S.-M. Liu and G.~Papavassilopoulos, ``Numerical experience with parallel
  algorithms for solving the {BMI} problem,'' \emph{IFAC-PapersOnline},
  vol.~29, no.~1, pp. 1827--1832, 1996.

\bibitem{hassibi1999path}
A.~Hassibi, J.~How, and S.~Boyd, ``A path-following method for solving {BMI}
  problems in control,'' in \emph{Proc. IEEE Am. Control Conf.}, vol.~2, 1999,
  pp. 1385--1389.

\bibitem{ibaraki2001rank}
S.~Ibaraki and M.~Tomizuka, ``Rank minimization approach for solving {BMI}
  problems with random search,'' in \emph{Proc. IEEE Am. Control Conf.},
  vol.~3, 2001, pp. 1870--1875.

\bibitem{doelman2016sequential}
R.~Doelman and M.~Verhaegen, ``Sequential convex relaxation for convex
  optimization with bilinear matrix equalities,'' in \emph{Proc. IEEE P.
  Control Conf.}, 2016, pp. 1946--1951.

\bibitem{beck2010sequential}
A.~Beck, A.~Ben-Tal, and L.~Tetruashvili, ``A sequential parametric convex
  approximation method with applications to nonconvex truss topology design
  problems,'' \emph{J. Global Optim.}, vol.~47, no.~1, pp. 29--51, 2010.

\bibitem{lee2016sequential}
D.~Lee and J.~Hu, ``A sequential parametric convex approximation method for
  solving bilinear matrix inequalities,'' in \emph{Proc. IEEE Conf. Decis.
  Control}, 2016, pp. 1965--1970.

\bibitem{dinh2012combining}
Q.~T. Dinh, S.~Gumussoy, W.~Michiels, and M.~Diehl, ``Combining convex--concave
  decompositions and linearization approaches for solving {BMI}s, with
  application to static output feedback,'' \emph{IEEE Trans. Autom. Control},
  vol.~57, no.~6, pp. 1377--1390, 2012.

\bibitem{el1997cone}
L.~El~Ghaoui, F.~Oustry, and M.~AitRami, ``A cone complementarity linearization
  algorithm for static output-feedback and related problems,'' \emph{IEEE
  Trans. Autom. Control}, vol.~42, no.~8, pp. 1171--1176, 1997.

\bibitem{fukuda2001branch}
M.~Fukuda and M.~Kojima, ``Branch-and-cut algorithms for the bilinear matrix
  inequality {E}igenvalue problem,'' \emph{Comput. Optim. Appl.}, vol.~19,
  no.~1, pp. 79--105, 2001.

\bibitem{tuan1999new}
D.~Tuan, P.~Apkarian, and Y.~Nakashima, ``A new {L}agrangian dual global
  optimization algorithm for solving bilinear matrix inequalities,'' in
  \emph{Proc. IEEE Am. Control Conf.}, vol.~3, 1999, pp. 1851--1855.

\bibitem{tuan2000low}
H.~D. Tuan and P.~Apkarian, ``Low nonconvexity-rank bilinear matrix
  inequalities: algorithms and applications in robust controller and structure
  designs,'' \emph{IEEE Trans. Autom. Control}, vol.~45, no.~11, pp.
  2111--2117, 2000.

\bibitem{chiu2017method}
W.-Y. Chiu, ``Method of reduction of variables for bilinear matrix inequality
  problems in system and control designs,'' \emph{IEEE Trans. Syst., Man,
  Cybern., Syst.}, vol.~47, no.~7, pp. 1241--1256, 2017.

\bibitem{kim2005generalized}
S.~Kim, M.~Kojima, and H.~Waki, ``Generalized {L}agrangian duals and sums of
  squares relaxations of sparse polynomial optimization problems,'' \emph{SIAM
  J. Optim.}, vol.~15, no.~3, pp. 697--719, 2005.

\bibitem{ahmadi2015sum}
A.~A. Ahmadi and G.~Hall, ``Sum of squares basis pursuit with linear and second
  order cone programming,'' \emph{Contemp. Math.}, 2017.

\bibitem{lasserre2001global}
J.~B. Lasserre, ``Global optimization with polynomials and the problem of
  moments,'' \emph{SIAM J. Optim.}, vol.~11, no.~3, pp. 796--817, 2001.

\bibitem{prajna2004sostools}
S.~Prajna, A.~Papachristodoulou, P.~Seiler, and P.~Parrilo, ``{SOSTOOLS}: Sum
  of squares optimization toolbox for {MATLAB}. user’s guide, version 2.00,''
  2004.

\bibitem{papachristodoulou2013sostools}
A.~Papachristodoulou, J.~Anderson, G.~Valmorbida, S.~Prajna, P.~Seiler, and
  P.~Parrilo, \emph{{SOSTOOLS}: Sum of squares optimization toolbox for
  {MATLAB}}, \texttt{http://arxiv.org/abs/1310.4716}, 2013, available from
  \texttt{http://www.eng.ox.ac.uk/control/sostools},
  \texttt{http://www.cds.caltech.edu/sostools} and
  \texttt{http://www.mit.edu/\~{}parrilo/sostools}.

\bibitem{henrion2003gloptipoly}
D.~Henrion and J.-B. Lasserre, ``{GloptiPoly}: Global optimization over
  polynomials with {MATLAB} and {SeDuMi},'' \emph{ACM Trans. Math. Softw.},
  vol.~29, no.~2, pp. 165--194, 2003.

\bibitem{shapiro1997first}
A.~Shapiro, ``First and second order analysis of nonlinear semidefinite
  programs,'' \emph{Math. Program.}, vol.~77, no.~1, pp. 301--320, 1997.

\end{thebibliography}
\appendix
In order to prove \cref{thm:BMI_thm_01,thm:BMI_thm_02}, we need to consider the following non-convex optimization problem
\begin{subequations} \label{eq:BMI_BMIOrig_Pen}
\begin{align}
&\underset{{\xbf}\in{\Rbb}^{n}}{\text{minimize}} && \cbf^{\!\top}\xbf+\eta \lVert\xbf-\check{\xbf}{\rVert}_{2}^2 \label{eq:BMI_BMIOrig_Pen_obj}\\
&\text{subject to}
&&p(\xbf,\xbf\xbf^{\!\top})\preceq 0, \label{eq:BMI_BMIOrig_Pen_con_01}
\end{align}
\end{subequations}
where $\check{\xbf}\in\Rbb^n$ is the initial point. Observe that problems \cref{eq:BMI_BMIOrig_obj,eq:BMI_BMIOrig_con_01} and \cref{eq:BMI_BMIOrig_Pen_obj,eq:BMI_BMIOrig_Pen_con_01} have the same feasible set, which is denoted by $\Fcal$. Assume that $\Fcal$ is nonempty with an arbitrary member $\xbf^{\prime}$. We define 
\begin{equation}
\begin{aligned}[b]
\!\!\Acal\!\triangleq\!\Big\{\xbf\!\in\!\Rbb^{n}\Big|\cbf^{\!\top}\!\xbf\msh+\msh\eta\lVert\xbf\!-\!\check{\xbf}{\rVert}_{2}^2\leq\cbf^{\!\top}\!\xbf^{\prime}\msh+\msh\eta\lVert\xbf^{\prime}\!-\!\check{\xbf}{\rVert}_{2}^2\Big\}.
\end{aligned}
\end{equation}
Due to the compactness of the set $\Acal\cap\Fcal$, it is straightforward to verify that the optimal solution of the problem \cref{eq:BMI_BMIOrig_Pen_obj,eq:BMI_BMIOrig_Pen_con_01} is attainable if $\eta>0$. 

\vspace{1.5mm}
\begin{lemma}\label{lm:BMI_Lemma_01}
Given an arbitrary $\varepsilon>0$, every optimal solution $\accentset{\ast}{\xbf}$ of problem \cref{eq:BMI_BMIOrig_Pen_obj,eq:BMI_BMIOrig_Pen_con_01} satisfies
\begin{equation}\label{eq:BMI_Lemma_01}
\begin{aligned}[b]
0\leq\lVert\accentset{\ast}{\xbf}-\check{\xbf}{\rVert}_{2}-d_{\sm{\Fcal}{5.5}}(\check{\xbf})\leq\varepsilon,
\end{aligned}
\end{equation}
if $\eta$ is sufficiently large. 
\end{lemma}
\vspace{1.5mm}

\begin{proof}
Consider an optimal solution $\accentset{\ast}{\xbf}$. Due to the \cref{def:DistFromSet}, the distance between $\check{\xbf}$ and every member of $\Fcal$ is greater than or equal to $d_{\sm{\Fcal}{5.5}}(\check{\xbf})$. Hence, the following inequality holds:
\begin{equation}\label{eq:BMI_DistRelLow}
\begin{aligned}[b]
0\leq\lVert\accentset{\ast}{\xbf}-\check{\xbf}{\rVert}_{2}-d_{\sm{\Fcal}{5.5}}(\check{\xbf}).
\end{aligned}
\end{equation}
Let $\xbf_h$ be an arbitrary member of $\{\xbf\msh\in\Fcal\vert\lVert{\xbf}-\check{\xbf}{\rVert}_{2}\msh=\msh d_{\sm{\Fcal}{5.5}}(\check{\xbf})\}$.
Due to the optimality of $\accentset{\ast}{\xbf}$, we have:
\begin{align}
\cbf^{\!\top}\accentset{\ast}{\xbf}+\eta\lVert\accentset{\ast}{\xbf}-\check{\xbf}{\rVert}_{2}^2&\leq\cbf^{\!\top}{\xbf}_{\sm{h}{6}}+\eta\lVert{\xbf}_{\sm{h}{6}}-\check{\xbf}{\rVert}_{2}^2, \label{eq:BMI_DistRelUp_01}
\end{align}
which implies that
\begin{align}
\Big\lVert(\accentset{\ast}{\xbf}-\check{\xbf})+\frac{\cbf}{2\eta}{\Big\rVert}_{2}&\leq\Big\lVert({\xbf}_{\sm{h}{6}}-\check{\xbf})+\frac{\cbf}{2\eta}{\Big\rVert}_{2}.\label{eq:BMI_DistRelUp_03}
\end{align}
Using the triangle inequality, we have 
\begin{align}
{\left\|\accentset{\ast}{\xbf}-\check{\xbf}\right\|}_{2}-\frac{1}{2\eta}\lVert\cbf{\rVert}_{2}&\leq\lVert{\xbf}_{\sm{h}{6}}\!-\!\check{\xbf}{\rVert}_{2}+\frac{1}{2\eta}\lVert\cbf{\rVert}_{2},\label{eq:BMI_DistBound_01} 
\end{align}
which leads to the following upper-bound:
\begin{equation}\label{eq:BMI_FinBound}
\begin{aligned}[b]
\lVert\accentset{\ast}{\xbf}-\check{\xbf}{\rVert}_{2}-d_{\sm{\Fcal}{5.5}}(\check{\xbf})\leq\frac{1}{\eta}\lVert\cbf{\rVert}_{2}.
\end{aligned}
\end{equation}
Hence, if
$\eta\geq\frac{\lVert\cbf{\rVert}_2}{\varepsilon}$, the combination of \cref{eq:BMI_DistRelLow} and \cref{eq:BMI_FinBound} completes the proof.
\end{proof}
\vspace{1.5mm}
In what follows, we obtain sufficient conditions to ensure that every solution of \cref{eq:BMI_BMIOrig_Pen_obj,eq:BMI_BMIOrig_Pen_con_01} satisfy the MFCQ condition. 

\vspace{1.5mm}
\begin{lemma}\label{lm:BMI_Lemma_02}
Assume that $\check{\xbf}\in\Rbb^{n}$ is a feasible point for \cref{eq:BMI_BMIOrig_Pen_obj,eq:BMI_BMIOrig_Pen_con_01} that satisfies the MFCQ condition. If $\eta$ is sufficiently large, every optimal solution $\accentset{\ast}{\xbf}$ of \cref{eq:BMI_BMIOrig_Pen_obj,eq:BMI_BMIOrig_Pen_con_01}, satisfies the MFCQ condition as well.
\end{lemma}
\vspace{1.5mm}

\begin{proof}
Consider an optimal solution $\accentset{\ast}{\xbf}$. Since the MFCQ condition holds for $\check{\xbf}$, there exists $\check{\bbf}\in\Rbb^n$ for which the conic inequality $p(\check{\xbf},\check{\xbf}\check{\xbf}^{\!\top})+\sum_{\sm{k\in\Ncal}{5.5}} \check{b}_{k}(\Kbf_{k}\!+\!{\delta}_{k}(\check{\xbf}))\prec 0$ is satisfied. Hence, due to the continuity of the matrix pencil $p$, if $\varepsilon$ is sufficiently small in \cref{lm:BMI_Lemma_01}, we have 
\begin{align}
p(\accentset{\ast}{\xbf},\accentset{\ast}{\xbf}\accentset{\ast}{\xbf}^{\!\top})+\sum_{\sm{k\in\Ncal}{5.5}} {\accentset{\ast}{b}}_{k}(\Kbf_{k}\!+\!{\delta}_{k}(\accentset{\ast}{\xbf}))\prec 0
\end{align}
which concludes the MFCQ condition holds for $\accentset{\ast}{\xbf}$.
\end{proof}
\vspace{1.5mm}

\begin{definition}
	Given an arbitrary symmetric matrix $\Lambdabf\in\Sbb_m$, define the matrix function $\alpha:\Sbb_m\rightarrow\Sbb_n$ as,
	\begin{equation}\label{eq:BMI_DefAlpha}
	\begin{aligned}[b]
	\alpha(\Lambdabf)&\triangleq\big[\langle\Lbf_{ij},\Lambdabf\rangle\big]_{\sm{ij\in\Ncal^2}{5.5}}.
	\end{aligned}
	\end{equation}
\end{definition}
\vspace{1.5mm}

It is straightforward to verify that
\begin{equation}\label{eq:BMI_DefAlphaProp}
\begin{aligned}[b]
2\alpha(\Lambdabf)\xbf=\sum_{\sm{k\in\Ncal}{6}}\langle\delta_{k}(\xbf),\Lambdabf\rangle{\ebf}_{k},
\end{aligned}
\end{equation}
for every $\xbf\in\Rbb^n$. This property will be used later in this section.
\vspace{1.5mm}

\vspace{1.5mm}
\begin{lemma}\label{lm:BMI_Lemma_03}
Assume that $\check{\xbf}\in\Rbb^n$ satisfies
\begin{align}\label{reqs}
s(\check{\xbf})> 2\lVert p{\rVert}_2 d_{\sm{\Fcal}{5.5}}(\check{\xbf}).
\end{align}
Given an arbitrary $\varepsilon>0$, every optimal solution $\accentset{\ast}{\xbf}$ of the problem \cref{eq:BMI_BMIOrig_Pen_obj,eq:BMI_BMIOrig_Pen_con_01} satisfies the inequality
\begin{equation}
\begin{aligned}
s(\check{\xbf})-s(\accentset{\ast}{\xbf})\leq 2\lVert p{\rVert}_2 d_{\sm{\Fcal}{5.5}}(\check{\xbf}) + \varepsilon,
\end{aligned}
\end{equation}
as well as the MFCQ condition, if $\eta$ is sufficiently large. 
\end{lemma}
\vspace{1.5mm}

\begin{proof}
Due to the definition of $s$, there exists $\check{\bbf}\in\Rbb^n$ such that $\|\check{\bbf}\|_2=1$ and
\begin{align}
s(\check{\xbf})=\ubar{\lambda}\Big(\!\!-\!\!\sum_{\sm{k\in\Ncal}{5.5}}\check{b}_{k}(\Kbf_{k}\!+\!\delta_{k}(\check{\xbf}))\!\Big).
\end{align}
As a result,
\begin{subequations}
\begin{align}
\hspace{-0.45cm}s(\accentset{\ast}{\xbf})& \geq\ubar{\lambda}\Big(\!\!-\!\!\sum_{\sm{k\in\Ncal}{5.5}}\check{b}_{k}(\Kbf_{k}\!+\!\delta_{k}(\accentset{\ast}{\xbf}))\!\Big) \\
&=\ubar{\lambda}\Big(\!\!-\!\!\sum_{\sm{k\in\Ncal}{5.5}}\check{b}_{k}(\Kbf_{k}\!+\!\delta_{k}(\check{\xbf}))\!-\!\!\sum_{\sm{k\in\Ncal}{5.5}}\check{b}_{k}\delta_{k}(\accentset{\ast}{\xbf}\!-\!\check{\xbf})\!\Big) \\
&\geq s(\check{\xbf})-\Big\Vert\sum_{\sm{k\in\Ncal}{5.5}}\check{b}_{k}\delta_{k}(\accentset{\ast}{\xbf}-\check{\xbf}){\Big\Vert}_2.
\end{align}
\end{subequations}
Let $\ubf$ be the eigenvector corresponding to the largest eigenvalue of $-\sum_{\sm{k\in\Ncal}{5.5}}\check{b}_{k}\delta_{k}(\accentset{\ast}{\xbf}-\check{\xbf})$. Then,   
\begin{subequations}
\begin{align}
s(\check{\xbf})-s(\accentset{\ast}{\xbf})&\leq \Big\Vert\sum_{\sm{k\in\Ncal}{5.5}}\check{b}_{k}\delta_{k}(\accentset{\ast}{\xbf}-\check{\xbf}){\Big\Vert}_2\\
&=\Big\vert{ \ubf^{\!\top}\Big(\sum_{\sm{k\in\Ncal}{5.5}}\check{b}_{k}\delta_{k}(\accentset{\ast}{\xbf}-\check{\xbf}) \Big)\ubf}\Big\vert \\
&=\Big\vert{\sum_{\sm{k\in\Ncal}{5.5}}\check{b}_{k}\langle\delta_{k}(\accentset{\ast}{\xbf}-\check{\xbf}),\ubf\ubf^{\!\top}\rangle}\Big\vert \\
&=\big\vert{\check{\bbf}^{\top}[\langle\delta_{k}(\accentset{\ast}{\xbf}-\check{\xbf}),\ubf\ubf^{\!\top}\rangle]_{\sm{k\in\Ncal}{5.5}}}\big\vert\\
&\leq\big\|[\langle\delta_{k}(\accentset{\ast}{\xbf}-\check{\xbf}),\ubf\ubf^{\!\top}\rangle]_{\sm{k\in\Ncal}{5.5}}\big\|_2.
\end{align}
\end{subequations}
On the other hand, according to the equation \eqref{eq:BMI_DefAlphaProp}, we have
\begin{align}
[\langle\delta_{k}(\accentset{\ast}{\xbf}-\check{\xbf}),\ubf\ubf^{\!\top}\rangle]_{\sm{k\in\Ncal}{5.5}}&=
2\alpha(\ubf\ubf^{\!\top})(\accentset{\ast}{\xbf}-\check{\xbf}),
\end{align}
which implies that
\begin{subequations}
\begin{align}
s(\check{\xbf})-s(\accentset{\ast}{\xbf})&\leq 2 \lVert\alpha(\ubf\ubf^{\!\top}){\rVert}_{2} \lVert\accentset{\ast}{\xbf}-\check{\xbf}{\rVert}_{2}\\
&\leq 2 \lVert p \rVert_{2} \lVert\accentset{\ast}{\xbf}-\check{\xbf}{\rVert}_{2}.
\end{align}
\end{subequations}
Therefore, according to \cref{lm:BMI_Lemma_01}, we have
\begin{align}
s(\check{\xbf})-s(\accentset{\ast}{\xbf})\leq 2 \lVert p \rVert_{2} \lVert\accentset{\ast}{\xbf}-\check{\xbf}{\rVert}_{2}
\leq 2 \lVert p \rVert_{2} d_{\sm{\Fcal}{5.5}}(\check{\xbf})+\varepsilon,\nonumber
\end{align}
if $\eta$ is sufficiently large. 

Additionally, for sufficiently small choices of $\varepsilon$, we have $s(\accentset{\ast}{\xbf})>0$. Hence, there exists $\accentset{\ast}{\bbf}\in\Rbb^n$ such that $\sum_{\sm{k\in\Ncal}{5.5}}\accentset{\ast}{b}_{k}(\Kbf_{k}\!+\!\delta_{k}(\accentset{\ast}{\xbf}))\prec 0$ and due to the feasibility of $\accentset{\ast}{\xbf}$, we have:
\begin{align}
p(\accentset{\ast}{\xbf},\accentset{\ast}{\xbf}\accentset{\ast}{\xbf}^\top)+\sum_{\sm{k\in\Ncal}{5.5}}\accentset{\ast}{b}_{k}(\Kbf_{k}\!+\!\delta_{k}(\accentset{\ast}{\xbf}))\prec 0,
\end{align}
which concludes the MFCQ condition for $\accentset{\ast}{\xbf}$.
\end{proof}
\vspace{1.5mm}
The following lemma ensures the existence of a dual certificate matrix, if the optimal solution of \cref{eq:BMI_BMIOrig_Pen_obj,eq:BMI_BMIOrig_Pen_con_01} satisfies the MFCQ condition.
\vspace{1.5mm}
\begin{lemma}\label{lm:BMI_Lemma_04}
	For every optimal solution $\accentset{\ast}{\xbf}$ of \cref{eq:BMI_BMIOrig_Pen_obj,eq:BMI_BMIOrig_Pen_con_01} which meets the MFCQ condition, there exists a dual matrix $\accentset{\ast}{\Lambdabf}\succeq 0$ such that the point $(\accentset{\ast}{\xbf},\accentset{\ast}{\Lambdabf})$ satisfies the following Karush-Kuhn-Tucker (KKT) equations:
	\begin{subequations}
		\begin{align}
			\cbf\msh+\msh2\eta(\accentset{\ast}{\xbf}-\check{\xbf})\msh+\msh\sum_{\sm{k\in\Ncal}{5.5}}\langle\Kbf_{k},\accentset{\ast}{\Lambdabf}\rangle\ebf_{k}\msh+\msh2\alpha(\accentset{\ast}{\Lambdabf})\accentset{\ast}{\xbf} &=0, \label{eq:BMI_KKTCond_01} \\
			\accentset{\ast}{\Lambdabf}p(\accentset{\ast}{\xbf},\accentset{\ast}{\xbf}{\accentset{\ast}{\xbf}}^{\!\top})&=0. \label{eq:BMI_KKTCond_02}
		\end{align}
	\end{subequations}
\end{lemma}
\vspace{1.5mm}

\begin{proof}
	Since the optimal solution $\accentset{\ast}{\xbf}$ satisfies the MFCQ condition, there exists a dual matrix $\accentset{\ast}{\Lambdabf}\succeq 0$ such that the point $(\accentset{\ast}{\xbf},\accentset{\ast}{\Lambdabf})$ satisfies the following conditions:
	\begin{subequations}
		\begin{align}
			\nabla_{\!\xbf}\Lcal_{\sm{\mathrm{p}}{6}}(\accentset{\ast}{\xbf},\accentset{\ast}{\Lambdabf})&=0, \label{eq:BMI_MFCond_01} \\
			\accentset{\ast}{\Lambdabf} p(\accentset{\ast}{\xbf},\accentset{\ast}{\xbf}{\accentset{\ast}{\xbf}}^{\!\top})&=0,\label{eq:BMI_MFCond_02}
		\end{align}
	\end{subequations}
	where $\nabla_{\!\xbf}$ represents the gradients with respect to $\xbf$ and $\Lcal_{\sm{\mathrm{p}}{6}}(\xbf,\Lambdabf)$ denotes the Lagrangian function of \cref{eq:BMI_BMIOrig_Pen_obj,eq:BMI_BMIOrig_Pen_con_01}:
	\begin{equation}
	\begin{aligned}[b]
	\Lcal_{\sm{\mathrm{p}}{6}}(\xbf,\Lambdabf)=\cbf^{\!\top}\xbf+\eta \lVert\xbf-\check{\xbf}{\rVert}_{2}^2+\langle p(\xbf,{\xbf}{\xbf}^{\!\top}),\Lambdabf\rangle.
	\end{aligned}
	\end{equation}
	Observe that \cref{eq:BMI_KKTCond_01,eq:BMI_KKTCond_02} and \cref{eq:BMI_MFCond_01,eq:BMI_MFCond_02} are equivalent. Therefore, the point $(\accentset{\ast}{\xbf}, \accentset{\ast}{\Lambdabf})$ satisfies the KKT conditions \cref{eq:BMI_KKTCond_01,eq:BMI_KKTCond_02}. 
\end{proof}
\vspace{1.5mm}
In following two lemmas bound the value of $\frac{\tr\{\accentset{\ast}{\Lambdabf}\}}{\eta}$ for both cases where $\check{\xbf}$ is feasible and infeasible.
\vspace{1.5mm}
\begin{lemma}\label{lm:BMI_Lemma_05}
Consider an arbitrary $\varepsilon>0$ and assume that $\check{\xbf}\in\Fcal$ is a feasible point for \cref{eq:BMI_BMIOrig_Pen_obj,eq:BMI_BMIOrig_Pen_con_01} that satisfies the MFCQ condition. If $\eta$ is sufficiently large, for every optimal solution $\accentset{\ast}{\xbf}$ of \cref{eq:BMI_BMIOrig_Pen_obj,eq:BMI_BMIOrig_Pen_con_01}, there exists a dual matrix $\accentset{\ast}{\Lambdabf}\succeq 0$ that satisfies the inequality
\begin{equation}\label{eq:BMI_Lemma_05_a}
\begin{aligned}
\frac{\tr\{\accentset{\ast}{\Lambdabf}\}}{\eta}\leq\varepsilon,
\end{aligned}
\end{equation}
as well as the equations \cref{eq:BMI_KKTCond_01,eq:BMI_KKTCond_02}.
\end{lemma}

\vspace{1.5mm}

\begin{proof}
According to \cref{lm:BMI_Lemma_02}, if $\eta$ is large enough, $\accentset{\ast}{\xbf}$ satisfies the MFCQ condition. Hence, there exists $\accentset{\ast}{\bbf}\in\Rbb^n$ such that 
\begin{align}\label{bbbb}
-p(\accentset{\ast}{\xbf},\accentset{\ast}{\xbf}{\accentset{\ast}{\xbf}}^{\!\top})-\!\sum_{\sm{k\in\Ncal}{5.5}}\accentset{\ast}{b}_{k}(\Kbf_{k}+\delta_{k}(\accentset{\ast}{\xbf}))\succeq0.
\end{align}
In addition, according to \cref{lm:BMI_Lemma_04} there exists $\accentset{\ast}{\Lambdabf}\succeq0$ such that the pair $(\accentset{\ast}{\xbf},\accentset{\ast}{\Lambdabf})$ satisfies the KKT equations \cref{eq:BMI_KKTCond_01,eq:BMI_KKTCond_02}. Therefore, pre-multiplying $\accentset{\ast}{\bbf}^{\top}$ to both sides of \cref{eq:BMI_KKTCond_01} yields:
\begin{equation}\label{aaaa}
\begin{aligned}[b]
\accentset{\ast}{\bbf}^{\top}(\cbf\!+\!2\eta(\accentset{\ast}{\xbf}\!-\!\check{\xbf})) \!+\!\langle\sum_{\sm{k\in\Ncal}{5.5}}\accentset{\ast}{b}_{k}(\Kbf_{k}\!+\!\delta_{k}(\accentset{\ast}{\xbf})),\accentset{\ast}{\Lambdabf}\rangle\!=\!0.
\end{aligned}
\end{equation}
Due to the matrix inequality \eqref{bbbb} and since $\accentset{\ast}{\Lambdabf}\succeq0$, we have: 
\begin{align}\label{nnnn}
\langle-p(\accentset{\ast}{\xbf},\accentset{\ast}{\xbf}{\accentset{\ast}{\xbf}}^{\!\top})-\!\sum_{\sm{k\in\Ncal}{5.5}}\accentset{\ast}{b}_{k}(\Kbf_{k}+\delta_{k}(\accentset{\ast}{\xbf})),\accentset{\ast}{\Lambdabf}\rangle\geq 0.
\end{align}
Hence, according to the complementary slackness \cref{eq:BMI_KKTCond_02}, we have: 
\begin{subequations}
\begin{align}
\!\!\!\!\tr\{\accentset{\ast}{\Lambdabf}\}s(\accentset{\ast}{\xbf})&\leq\langle-p(\accentset{\ast}{\xbf},\accentset{\ast}{\xbf}{\accentset{\ast}{\xbf}}^{\!\top})-\!\sum_{\sm{k\in\Ncal}{5.5}}\accentset{\ast}{b}_{k}(\Kbf_{k}+\delta_{k}(\accentset{\ast}{\xbf})),\accentset{\ast}{\Lambdabf}\rangle \label{eq:BMI_PRO_LM04_00} \\
&=\langle-\!\sum_{\sm{k\in\Ncal}{5.5}}\accentset{\ast}{b}_{k}(\Kbf_{k}+\delta_{k}(\accentset{\ast}{\xbf})),\accentset{\ast}{\Lambdabf}\rangle \label{eq:BMI_PRO_LM04_01} \\
&=\langle\accentset{\ast}{\bbf},\cbf+2\eta(\accentset{\ast}{\xbf}-\check{\xbf})\rangle \label{eq:BMI_PRO_LM04_02}\\
&=\langle\accentset{\ast}{\bbf},\cbf\rangle+2\eta\langle\accentset{\ast}{\bbf},(\accentset{\ast}{\xbf}-\check{\xbf})\rangle\label{eq:BMI_PRO_LM04_03}\\
&\leq\lVert\accentset{\ast}{\bbf}{\rVert}_{2}\lVert\cbf{\rVert}_{2}+2\eta\lVert\accentset{\ast}{\bbf}{\rVert}_{2}\lVert\accentset{\ast}{\xbf}-\check{\xbf}{\rVert}_{2}\label{eq:BMI_PRO_LM04_04}\\
&=\lVert\cbf{\rVert}_{2}+2\eta\lVert\accentset{\ast}{\xbf}-\check{\xbf}{\rVert}_{2},\label{eq:BMI_PRO_LM04_05}
\end{align}
\end{subequations}
and therefore:
\begin{align}\label{ineq}
\frac{\tr\{\accentset{\ast}{\Lambdabf}\}}{\eta}&\leq\frac{\lVert\cbf{\rVert}_{2}}{\eta s(\accentset{\ast}{\xbf})}+
\frac{2\lVert\accentset{\ast}{\xbf}-\check{\xbf}{\rVert}_{2}}{s(\accentset{\ast}{\xbf})}.
\end{align}
According to \cref{lm:BMI_Lemma_01}, if $\eta$ is large,
$\lVert\accentset{\ast}{\xbf}-\check{\xbf}{\rVert}_{2}$ is arbitrarily small. Due to the continuity of $s$, we can argue that 
$| s(\accentset{\ast}{\xbf})-s(\check{\xbf})|$ is arbitrarily small as well. Now, since $s(\check{\xbf})>0$, the right side of the inequality \eqref{ineq} is not greater than $\varepsilon$, if $\eta$ is sufficiently large.
\end{proof}
\vspace{1.5mm}

\vspace{1.5mm}
\begin{lemma}\label{lm:BMI_Lemma_06}
Consider an arbitrary $\varepsilon>0$ and assume that $\check{\xbf}\in\Rbb^n$ satisfies the inequality \eqref{reqs}.
If $\eta$ is sufficiently large, for every optimal solution $\accentset{\ast}{\xbf}$ of \cref{eq:BMI_BMIOrig_Pen_obj,eq:BMI_BMIOrig_Pen_con_01}, there exists a dual matrix $\accentset{\ast}{\Lambdabf}\succeq 0$ that satisfies the inequality	
\begin{equation}\label{eqkey}
\begin{aligned}
\frac{\tr\{\accentset{\ast}{\Lambdabf}\}}{\eta} &\leq\frac{2d_{\sm{\Fcal}{5.5}}(\check{\xbf})}{s(\check{\xbf})-2\lVert p{\rVert}_2d_{\sm{\Fcal}{5.5}}(\check{\xbf})}+\varepsilon,
\end{aligned}
\end{equation}
as well as the equations \cref{eq:BMI_KKTCond_01,eq:BMI_KKTCond_02}.
\end{lemma}

\vspace{1.5mm}
\begin{proof}
According to the \cref{lm:BMI_Lemma_03}, $\accentset{\ast}{\xbf}$ satisfies the MFCQ condition. In addition, the \cref{lm:BMI_Lemma_04} implies that there exists $\accentset{\ast}{\Lambdabf}\succeq0$ such that point $(\accentset{\ast}{\xbf},\accentset{\ast}{\Lambdabf})$ satisfies the KKT equations \cref{eq:BMI_KKTCond_01,eq:BMI_KKTCond_02}. Since $s(\accentset{\ast}{\xbf})>0$, we can similarly argue that the inequality \eqref{ineq} holds true:
\begin{align}\label{ineq2}
\!\!\frac{\tr\{\accentset{\ast}{\msh\Lambdabf\msh}\}}{\eta}\!\leq\!\frac{\lVert\cbf{\rVert}_{2}}{\eta s(\accentset{\ast}{\xbf})}\!+\!
	\frac{2\lVert\accentset{\ast}{\xbf}\msh-\msh\check{\xbf}{\rVert}_{2}}{s(\accentset{\ast}{\xbf})}\!\leq\!
		\frac{\lVert\cbf{\rVert}_{2}\msh+\msh2\eta\lVert\accentset{\ast}{\xbf}\msh-\msh\check{\xbf}{\rVert}_{2}}{\eta[s(\check{\xbf})\!-\!2\lVert p{\rVert}_2d_{\sm{\Fcal}{5.5}}(\check{\xbf})]}.
\end{align}
Now, according to \cref{lm:BMI_Lemma_01}, if $\eta$ is large, the above inequality concludes \eqref{eqkey}.
\end{proof}
\vspace{1.5mm}

The next lemma presents sufficient conditions under which the optimal solution of \cref{eq:BMI_BMIOrig_Pen_obj,eq:BMI_BMIOrig_Pen_con_01} can be obtained by solving penalized convex relaxation.
\vspace{1.5mm}
\begin{lemma}\label{lm:BMI_Lemma_07}
Consider an optimal solution $\accentset{\ast}{\xbf}\in\mathcal{F}$ for the problem \cref{eq:BMI_BMIOrig_Pen_obj,eq:BMI_BMIOrig_Pen_con_01}, and a matrix $\accentset{\ast}{\Lambdabf}\succeq 0$ such that point $(\accentset{\ast}{\xbf},\accentset{\ast}{\Lambdabf})$ satisfies the  conditions \cref{eq:BMI_KKTCond_01,eq:BMI_KKTCond_02}. Then, the pair $(\accentset{\ast}{\xbf},\accentset{\ast}{\xbf}{\accentset{\ast}{\xbf}}^{\!\top})$ is the unique primal solution to the penalized convex relaxation problem \cref{eq:BMI_GenRelaxation_Pen_obj,eq:BMI_GenRelaxation_Pen_con_01,eq:BMI_GenRelaxation_Pen_con_02}, if the following conic inequality holds true:
\begin{equation}\label{eq:BMI_DisCond}
\begin{aligned}
\eta\Ibf+\alpha(\accentset{\ast}{\Lambdabf})\succ_{\Ccal^{\ast}_{\msh\sm{k}{5}}} 0,
\end{aligned}
\end{equation}
where $k\in\{1,2,3\}$, $\Ccal^{\ast}_{\msh\sm{k}{5}}$ denotes the dual cone of $\Ccal_k$, and the matrices $\accentset{\ast}{\Lambdabf}$ and $\eta\Ibf+\alpha(\accentset{\ast}{\Lambdabf})$ are the dual optimal Lagrange multipliers associated with the constraints
\cref{eq:BMI_GenRelaxation_Pen_con_01} and \cref{eq:BMI_GenRelaxation_Pen_con_02}, respectively.
\end{lemma}

\vspace{1.5mm}

\begin{proof}
The Lagrangian of the penalized relaxation problem \cref{eq:BMI_GenRelaxation_Pen_obj,eq:BMI_GenRelaxation_Pen_con_01,eq:BMI_GenRelaxation_Pen_con_02} can be formed as follows,
\begin{align}
\Lcal_{\sm{\mathrm{r}}{6}}(\xbf,\Xbf,\Lambdabf)&=\cbf^{\!\top}\xbf+\eta\langle\Xbf-2\xbf\check{\xbf}^{\!\top},\Ibf\rangle+\langle p(\xbf,\Xbf),\Lambdabf\rangle\nonumber\\
&-\langle\eta\Ibf+\alpha(\accentset{\ast}{\Lambdabf}),\Xbf-\xbf\xbf^{\!\top}\rangle,
\end{align}
where $\eta\Ibf+\alpha(\accentset{\ast}{\Lambdabf})\in\Ccal^{\ast}_{\msh\sm{k}{5}}$ is the dual variable associated with \cref{eq:BMI_GenRelaxation_Pen_con_02}. 
Due to the convexity of the penalized relaxation problem, if a pair $\big((\accentset{\ast}{\xbf},\accentset{\ast}{\Xbf}),\accentset{\ast}{\Lambdabf}\big)$ satisfies the KKT conditions
\begin{subequations}
\begin{align}
\cbf-2\eta\check{\xbf}+\sum_{\sm{k\in\Ncal}{5.5}}\langle\Kbf_{k},\accentset{\ast}{\Lambdabf}\rangle\ebf_{k}+2(\eta\Ibf+\alpha(\accentset{\ast}{\Lambdabf}))\accentset{\ast}{\xbf}&=0,\label{eq:BMI_PenRel_KKTCond_01}\\[-0.75em]
\langle p(\accentset{\ast}{\xbf},\accentset{\ast}{\Xbf}),\accentset{\ast}{\Lambdabf}\rangle&=0,\label{eq:BMI_PenRel_KKTCond_02}\\
p(\accentset{\ast}{\xbf},\accentset{\ast}{\Xbf})&\preceq0,\label{eq:BMI_PenRel_KKTCond_03}\\
\eta\Ibf+\alpha(\accentset{\ast}{\Lambdabf})&\in\Ccal^{\ast}_{\msh\sm{k}{5}},\label{eq:BMI_PenRel_KKTCond_04}
\end{align}
\end{subequations}
then it is an optimal primal-dual solution for \cref{eq:BMI_GenRelaxation_Pen_obj,eq:BMI_GenRelaxation_Pen_con_01,eq:BMI_GenRelaxation_Pen_con_02}.

It can be easily verified that the KKT conditions
\cref{eq:BMI_PenRel_KKTCond_01,eq:BMI_PenRel_KKTCond_02,eq:BMI_PenRel_KKTCond_03,eq:BMI_PenRel_KKTCond_04} are satisfied for
$\big((\accentset{\ast}{\xbf},\accentset{\ast}{\xbf}\accentset{\ast}{\xbf}^{\top}),\accentset{\ast}{\Lambdabf}\big)$ as a direct consequence of  \cref{eq:BMI_KKTCond_01,eq:BMI_KKTCond_02}, \cref{eq:BMI_DisCond} and \cref{eq:BMI_BMIOrig_Pen_con_01}. Moreover, $(\accentset{\ast}{\xbf},\accentset{\ast}{\xbf}\accentset{\ast}{\xbf}^{\top})$ is the unique solution of the primal problem since $\eta\Ibf+\alpha(\accentset{\ast}{\Lambdabf})$ belongs to the interior of $\Ccal^{\ast}_{\msh\sm{k}{5}}$.
\end{proof}
\vspace{1.5mm}

\vspace{1.5mm}
\begin{lemma}\label{lm:BMI_Lemma_08}
Consider an optimal solution $\accentset{\ast}{\xbf}\in\mathcal{F}$ for problem \cref{eq:BMI_BMIOrig_Pen_obj,eq:BMI_BMIOrig_Pen_con_01}, and a matrix $\accentset{\ast}{\Lambdabf}\succeq 0$ such that point $(\accentset{\ast}{\xbf},\accentset{\ast}{\Lambdabf})$ satisfies the KKT equations \cref{eq:BMI_KKTCond_01,eq:BMI_KKTCond_02}. The pair $(\accentset{\ast}{\xbf},\accentset{\ast}{\xbf}{\accentset{\ast}{\xbf}}^{\!\top})$ is the unique primal solution to the penalized convex relaxation problem \cref{eq:BMI_GenRelaxation_Pen_obj,eq:BMI_GenRelaxation_Pen_con_01,eq:BMI_GenRelaxation_Pen_con_02}, if the following inequality holds true:
\begin{align}
\frac{\tr\{\accentset{\ast}{\Lambdabf}\}}{\eta}\leq\frac{\zeta_k}{\lVert p{\rVert}_2}
\end{align}
where $k\in\{1,2,3\}$, $\zeta_1=1$, $\zeta_2=(n-1)^{-1}$, and $\zeta_3=n^{-\frac{1}{2}}$.
\end{lemma}

\vspace{1.5mm}

\begin{proof}
According to \cref{lm:BMI_Lemma_05}, it suffices to verify \cref{eq:BMI_DisCond} in order to prove that $(\accentset{\ast}{\xbf},\accentset{\ast}{\xbf}{\accentset{\ast}{\xbf}}^{\!\top})$ is the unique optimal solution. Denote the eigenvalues and eigenvectors of $\accentset{\ast}{\Lambdabf}$ by $\{\accentset{\ast}{\lambda}_l\}_{\sm{l\in\Mcal}{5.5}}$ and $\{\accentset{\ast}{\ubf}_l\}_{\sm{l\in\Mcal}{5.5}}$, respectively. 
Hence:
\begin{subequations}\label{youoi}
	\begin{align}
	\big\lVert\alpha(\accentset{\ast}{\Lambdabf}){\big\rVert}_q&=\big\lVert\sum_{\sm{l\in\Mcal}{5.5}}\accentset{\ast}{\lambda}_l[\langle\Lbf_{ij},\accentset{\ast}{\ubf}_l\accentset{\ast}{\ubf}_l^{\!\top}{\rangle]}_{ij}{\big\rVert}_q\\
	&\leq\sum_{\sm{l\in\Mcal}{5.5}}\accentset{\ast}{\lambda}_l \big\lVert[\langle\Lbf_{ij},\accentset{\ast}{\ubf}_l\accentset{\ast}{\ubf}_l^{\!\top}{\rangle]}_{ij}{\big\rVert}_q\\
	&\leq\sum_{\sm{l\in\Mcal}{5.5}}\accentset{\ast}{\lambda}_l\big\lVert\alpha(\accentset{\ast}{\ubf}_l\accentset{\ast}{\ubf}_l^{\!\top}){\big\rVert}_q
	=\lVert p{\rVert}_q\tr\{\accentset{\ast}{\Lambdabf}\},
	\end{align}
\end{subequations}
\begin{enumerate}
\item SDP relaxation: The cone of positive semidefinite matrices is self-dual i.e., $\Ccal^{\ast}_{\sm{\mathrm{1}}{5}}=\Ccal_{\sm{\mathrm{1}}{5}}$. Therefore, in order to prove \cref{eq:BMI_DisCond}, it suffices to show that
\begin{align}
\eta-\big\lVert\alpha(\accentset{\ast}{\Lambdabf}){\big\rVert}_2\geq 0.
\end{align}
Hence, according to the bound provided in \eqref{youoi}, $(\accentset{\ast}{\xbf},\accentset{\ast}{\xbf}{\accentset{\ast}{\xbf}}^{\!\top})$ is the unique solution for the penalized SDP relaxation, if
\begin{equation}
\begin{aligned}[b]
\frac{\tr\{\accentset{\ast}{\Lambdabf}\}}{\eta}\leq\frac{1}{\lVert p{\rVert}_2}.
\end{aligned}
\end{equation}
\item SOCP relaxation: The dual cone $\Ccal_{\sm{\mathrm{2}}{5}}^{\ast}$ can be expressed as:
\begin{equation}
\begin{aligned}[b]
\hspace{-0.65cm}\Ccal_{\sm{\mathrm{2}}{5}}^{\ast}\!\triangleq&\Big\{\!\sum_{\begin{subarray}{c}\sm{i,j\in\Ncal}{5.5}\end{subarray}}\![\ebf_{i},\!\ebf_{j}]\,\Hbf_{ij}[\ebf_{i},\!\ebf_{j}]^{\!\top}\Big| \Hbf_{ij}\!\in\!\Sbb_{2}^{+},\sm{\forall i,\!j\!\in\!\Ncal}{8} \Big\}.\!\!\!\!\!\!\!\!\!\!\!\!
\end{aligned}
\end{equation}

Consider the following decomposition:
\begin{align}
\eta\Ibf+\alpha(\accentset{\ast}{\Lambdabf})=\sum_{\begin{subarray}{c}\sm{i,j\in\Ncal}{5.5}\\\sm{i\neq j}{5.5}\end{subarray}}\![\ebf_{i},\!\ebf_{j}]\,
\Abf_{ij}
\,[\ebf_{i},\!\ebf_{j}]^{\!\top},
\end{align}
where for every $(i,j)\in\Ncal^2$ we have
\begin{align}\small\!\!\!\!
\Abf_{ij}\!\triangleq\!\!
\begin{bmatrix}
\frac{\eta-[\alpha(\msh\accentset{\ast}{\Lambdabf}\msh)]_{ii}}{n-1} & \!\!\!-[\alpha(\msh\accentset{\ast}{\Lambdabf}\msh)]_{ij}\\
-[\alpha(\msh\accentset{\ast}{\Lambdabf}\msh)]_{ji} & \!\!\!\frac{\eta-[\alpha(\msh\accentset{\ast}{\Lambdabf}\msh)]_{jj}}{n-1}
\end{bmatrix}\!\msh \succeq \!\msh\left(\!\msh\frac{\eta}{n\!-\!1}\!-\! \|\alpha(\accentset{\ast}{\Lambdabf})\|_2\!\!\right)\!\Ibf_2\!\!\nonumber
\end{align}
Therefore, the inequality \cref{eq:BMI_DisCond} is satisfied for $\Ccal_{\sm{\mathrm{2}}{5}}^{\ast}$ if $\frac{\eta}{n-1}\geq \|\alpha(\accentset{\ast}{\Lambdabf})\|_2$. Now, according to the bound provided in \eqref{youoi}, $(\accentset{\ast}{\xbf},\accentset{\ast}{\xbf}{\accentset{\ast}{\xbf}}^{\!\top})$ is the unique optimal solution for the penalized SOCP relaxation if  
\begin{equation}
\begin{aligned}[b]
\frac{\tr\{\accentset{\ast}{\Lambdabf}\}}{\eta}\leq\frac{1}{(n-1)\lVert p{\rVert}_2}.
\end{aligned}
\end{equation}

\item Parabolic relaxation: The dual cone of $\Ccal_{\sm{\mathrm{3}}{5}}$ is the set of $n\times n$ symmetric diagonally dominant matrices defined as:
\begin{equation}
\Ccal_{\sm{\mathrm{3}}{5}}^{\ast} = \Big\{ \Hbf\in\Sbb_{n}\,\Big|\,|H_{ii}| \geq\!\!\!\!\sum_{\sm{j\in\Ncal\setminus\{i\}}{5.5}}\!\!|H_{ij}|,\,\sm{\forall i\!\in\!\Ncal}{8} \Big\}.
\end{equation}
Therefore, in order to prove \cref{eq:BMI_DisCond}, it suffices to show that
\begin{align}
\eta-\big\lVert\alpha(\accentset{\ast}{\Lambdabf}){\big\rVert}_1\geq 0.
\end{align}
Once again, the bound presented in \eqref{youoi} implies that $(\accentset{\ast}{\xbf},\accentset{\ast}{\xbf}{\accentset{\ast}{\xbf}}^{\!\top})$ is the unique solution for the penalized parabolic relaxation if $\frac{\tr\{\accentset{\ast}{\Lambdabf}\}}{\eta}\leq\frac{1}{\lVert p{\rVert}_1}$. This is a direct consequence of
\begin{align}
\frac{\tr\{\accentset{\ast}{\Lambdabf}\}}{\eta}\leq\frac{n^{-\frac{1}{2}}}{\lVert p{\rVert}_2}
\end{align}
since $\lVert p{\rVert}_1\leq n^{\frac{1}{2}}\lVert p{\rVert}_2$. 
\end{enumerate}
\end{proof}
\vspace{1.5mm}

\begin{proof}[\cref{thm:BMI_thm_01}]
Consider an arbitrary optimal solution $\accentset{\ast}{\xbf}$ for the problem \cref{eq:BMI_BMIOrig_Pen_obj,eq:BMI_BMIOrig_Pen_con_01}. According to \cref{lm:BMI_Lemma_04}, if $\eta$ is large enough, there exists a dual matrix $\accentset{\ast}{\Lambdabf}\succeq 0$ such that point $(\accentset{\ast}{\xbf},\accentset{\ast}{\Lambdabf})$ satisfies the KKT equations \cref{eq:BMI_KKTCond_01,eq:BMI_KKTCond_02}, as well as the inequality \cref{eq:BMI_Lemma_05_a} for $\varepsilon=\frac{\min\{\zeta_1,\zeta_2,\zeta_3\}}{2\|p\|_2}$. Therefore, according to the \cref{lm:BMI_Lemma_08}, the pair $(\accentset{\ast}{\xbf},\accentset{\ast}{\xbf}{\accentset{\ast}{\xbf}^{\!\top}}\!)$ is the unique primal solution 
to the penalized convex relaxation problem \cref{eq:BMI_GenRelaxation_Pen_obj,eq:BMI_GenRelaxation_Pen_con_01,eq:BMI_GenRelaxation_Pen_con_02}.
\end{proof}

\vspace{3mm}

\begin{proof}[\cref{thm:BMI_thm_02}]
Consider an arbitrary optimal solution $\accentset{\ast}{\xbf}$ for the problem \cref{eq:BMI_BMIOrig_Pen_obj,eq:BMI_BMIOrig_Pen_con_01}. According to the \cref{lm:BMI_Lemma_04}, if $\eta$ is sufficiently large, there exists a dual matrix $\accentset{\ast}{\Lambdabf}\succeq 0$ such that point $(\accentset{\ast}{\xbf},\accentset{\ast}{\Lambdabf})$ satisfies the KKT equations \cref{eq:BMI_KKTCond_01,eq:BMI_KKTCond_02}, as well as the inequality \cref{eqkey} for any arbitrarily $\varepsilon$.
It is straightforward to verify that
\begin{align}
\frac{d_{\sm{\Fcal}{5.5}}(\check{\xbf})}{s(\check{\xbf})} < \frac{\omega_k}{\lVert p{\rVert}_2}
\;\;\Rightarrow\;\;
\frac{2d_{\sm{\Fcal}{5.5}}(\check{\xbf})}{s(\check{\xbf})-2\lVert p{\rVert}_2d_{\sm{\Fcal}{5.5}}(\check{\xbf})}<\frac{\zeta_k}{\lVert p{\rVert}_2},
\end{align}
for all $k\in\{1,2,3\}$. Therefore, according to \cref{lm:BMI_Lemma_08}, the pair $(\accentset{\ast}{\xbf},\accentset{\ast}{\xbf}{\accentset{\ast}{\xbf}}^{\!\top})$ is the unique primal solution 
to the penalized convex relaxation problem \cref{eq:BMI_GenRelaxation_Pen_obj,eq:BMI_GenRelaxation_Pen_con_01,eq:BMI_GenRelaxation_Pen_con_02}.
\end{proof}

\end{document}